\documentclass[conference,9pt]{IEEEtran}

\usepackage{cite,url,graphicx,multirow}
\usepackage{amssymb, dsfont, bm}
\usepackage[cmex10]{amsmath}
\usepackage{empheq}
\usepackage{amsthm}

\newtheoremstyle{altieeeassumstyle}
  {5pt}
  {5pt}
  {}
  {}
  {\bfseries}
  {.}
  { }
  {\thmname{#1}\thmnumber{#2}\thmnote{ #3}}

\newtheoremstyle{altieeeremarkstyle}
  {5pt}
  {5pt}
  {}
  {}
  {\bfseries}
  {.}
  { }
  {\thmname{#1} \thmnumber{#2}\thmnote{ #3}}

\theoremstyle{altieeeassumstyle}

\theoremstyle{altieeeremarkstyle}
\newtheorem{remark}{Remark}

\newcommand{\cT}{\mathcal{T}}
\newcommand{\cN}{\mathcal{N}_a}
\newcommand{\cNg}{\mathcal{N}_g}

\DeclareMathOperator*{\argmin}{arg\,min}

\allowdisplaybreaks[4]

\IEEEoverridecommandlockouts

\begin{document}
\title{Disaggregated Bundle Methods for \\
Distributed Market Clearing in Power Networks}

\author{\IEEEauthorblockN{Yu Zhang\IEEEauthorrefmark{1},
Nikolaos Gatsis\IEEEauthorrefmark{2},
Georgios B. Giannakis\IEEEauthorrefmark{1}}
\IEEEauthorblockA{\IEEEauthorrefmark{1}Dept. of ECE and DTC,
University of Minnesota, Minneapolis, USA\\
Emails: \{zhan1220,georgios\}@umn.edu}
\IEEEauthorblockA{\IEEEauthorrefmark{2}Dept. of ECE,
The Univ. of Texas at San Antonio, San Antonio, USA \\
Email: Nikolaos.Gatsis@utsa.edu}
\thanks{This work was supported by NSF-ECCS grant 1202135,
and Institute of Renewable Energy and the Environment (IREE)
grant RL-0010-13, University of Minnesota.}}

\maketitle

\begin{abstract}
A fast distributed approach is developed for the market clearing with large-scale
demand response in electric power networks. In addition to conventional supply bids,
demand offers from aggregators serving large numbers of residential smart appliances with different
energy constraints are incorporated. Leveraging the Lagrangian relaxation based
dual decomposition, the resulting optimization problem is decomposed into
separate subproblems, and then solved in a distributed fashion by the market
operator and each aggregator aided by the end-user smart meters.
A disaggregated bundle method is adapted for solving the dual problem with a
separable structure. Compared with the conventional dual update algorithms,
the proposed approach exhibits faster convergence speed, which results in 
reduced communication overhead. Numerical results corroborate
the effectiveness of the novel approach.
\end{abstract}


\section{Introduction}
\label{sec:intro}

Demand response (DR) has been identified as an important resource management task
in modern power networks
promising to enable end-user interaction with the grid.
DR aggregators serving large numbers of residential users will be able to participate
in the market clearing by offering bids 
depending on the elasticity for power consumption of their end users.
%
%
%
Bidirectional communication between aggregators and users
is provided by the Advanced Metering Infrastructure (AMI)~\cite{wollen11acc},
with smart meters used as the end-users' terminals.

The principal challenge for large-scale incorporation of DR from residential end-users is
to account for the user scheduling preferences and intertemporal flexibility in a way that also protects user privacy.
The advantages of intertemporal load scheduling flexibility are for instance demonstrated in~\cite{SuKi09, WaKeKi10},
but without considering small-scale users and pertinent distributed algorithms.
Aggregation of small-scale user loads into the system scheduling has been the theme of~\cite{JoIl10a, TrPeBeSt11},
but an array of issues ranging from incorporation of user utility functions and user privacy to algorithm convergence,
are not fully addressed.
Algorithms for market clearing with large-scale integration of DR from small loads
with different utility functions are developed in~\cite{GaGG-TSG13} based on Lagrangian dual decomposition.
The disaggregated cutting plane method (CPM) is proposed therein for updating the Lagrange multipliers.

This paper proposes a market clearing approach distributed among the market operator, aggregators,
and the user smart meters by building upon the earlier work in~\cite{GaGG-TSG13}.
Each end-user has preferences for smart appliance scheduling captured by utility functions
and intertemporal constraints.
The objective is to minimize the social net cost for day-ahead
market clearing, while transmission network constraints are
included in the form of DC power flows.
To cope with the challenges of respecting end-user privacy and
large-scale DR, dual decomposition is applied to the
resulting optimization problem.
Leveraging Lagrangian relaxation of the coupling constraints,
the large-scale optimization decomposes into manageable small problems
solved by the market operator (MO) and the aggregators
in conjunction with the residential smart meters.
Exploiting the separable structure of the problem at hand,
a disaggregated bundle method is introduced for
solving the dual problem with guaranteed convergence of
the Lagrange multipliers.
The developed solver yields faster convergence than its CPM counterpart,
implying less communication overhead  between the MO and the aggregators.

The remainder of this paper is organized as follows.
Section~\ref{sec:MC} presents the market clearing problem involving large-scale DR.
The decomposition algorithm along with the disaggregated bundle method
solver is developed in Section~\ref{sec:decomp}.
Numerical tests are in Section~\ref{sec:num}, while conclusions and
future directions are offered in Section~\ref{sec:concl}.

\section{Market Clearing Formulation}\label{sec:MC}

Consider a power network comprising $N_g$ generators, $N_b$ buses, $N_l$ lines, and
$N_a$ aggregators, each serving a large number
of residential end-users with controllable smart appliances.
The scheduling horizon of interest is $\cT:=\{1,2,\ldots,T\}$
(e.g., one day ahead).
Let $\mathbf{p}_G^t := [P_{G_1}^t,\ldots,P_{G_{N_g}}^t]^{\prime}$
and $\mathbf{p}_{\mathrm{DRA}}^t :=
[P_{\mathrm{DRA}_1}^t,\ldots,P_{\mathrm{DRA}_{N_a}}^t]^{\prime}$
denote the generator power outputs, and the power consumption of
the aggregators at slot $t$, respectively.\footnote {$\mathbf{x}^{\prime}$ denotes transpose of the vector $\mathbf{x}$.}
Define further the sets $\cN:=\{1,2,\ldots,N_a\}$ and $\cNg:=\{1,2,\ldots,N_g\}$.
Each aggregator $j\in\cN$ serves a set $\mathcal{R}_j$ of residential users,
and each user $r\in\mathcal{R}_j$ has a set $\mathcal{S}_{rj}$ of controllable smart appliances.
Let $\mathbf{p}_{jrs}$ be the power consumption of smart appliance $s$
and user $r$ corresponding to aggregator $j$ across the horizon.
The power consumption $\mathbf{p}_{jrs}$ of each smart appliance across the horizon must typically satisfy
operating constraints captured by a set $\mathcal{P}_{jrs}$, and may also give rise to user satisfaction 
represented by a concave utility function $B_{jrs}(\mathbf{p}_{jrs})$.
Moreover, the generation cost is captured by convex functions $\{C_i(\cdot)\}_{i}$, and
the fixed base load demands across the network buses at slot $t$
is denoted by the vector $\mathbf{p}_{\mathrm{BL}}^t$.

For brevity, vector $\mathbf{p}_0$ is used to collect all $p_{G_i}^t$, $P_{\mathrm{DRA}_j}^t$, and network nodal angles $\theta_n^t$; while vector $\mathbf{p}_j$ ($j\in\cN$) collects all smart appliance consumptions corresponding to aggregator $j$.
With the goal of minimizing the system net cost, the DC optimal power flow (OPF)
based market clearing stands as follows:
\begin{subequations}
\label{eq:mc-ALL}
\begin{align}
f^{*}= &\min_{ \{\mathbf{p}_j\}_{j=0}^{N_a}}
\sum_{t=1}^T \sum_{i=1}^{N_g} C_{i}(P_{G_i}^t)
- \sum_{j=1}^{N_a} \sum_{r\in\mathcal{R}_{j}}
\sum_{s\in\mathcal{S}_{jr}} B_{jrs}(\mathbf{p}_{jrs}) \label{eq:mc-obj}\\
\text{s. t.} ~~&
\mathbf{A}_g \mathbf{p}_G^t -\mathbf{A}_a \mathbf{p}_{\mathrm{DRA}}^t
- \mathbf{p}_{\mathrm{BL}}^t = \mathbf{B} \bm\theta^t,~t \in \cT \label{eq:mc-bus} \\
& P_{G_i}^{\min} \leq P_{G_i}^t \leq P_{G_i}^{\max},
~i \in \cNg,~t \in \cT \label{eq:mc-genlim}\\
& -\mathsf{R}_i^{\mathrm{down}} \leq P_{G_i}^t - P_{G_i}^{t-1} \leq \mathsf{R}_i^{\mathrm{up}},
~i \in \cNg, \: t \in \cT \label{eq:mc-ramp}\\
& \mathbf{f}^{\min} \leq \mathbf{H} \bm\theta^t \leq \mathbf{f}^{\max},
~t \in \cT \label{eq:mc-flowlim}\\
& \theta_1^t=0, \: t \in \cT \label{eq:mc-refbus} \\
& 0 \leq P_{\mathrm{DRA}_j}^t \leq P_{\mathrm{DRA}_j}^{\max},
~j \in \cN,\,  t \in \cT \label{eq:mc-DRAlim}\\
& P_{\mathrm{DRA}_{j}}^t = \sum\nolimits_{r\in\mathcal{R}_{j},\,s\in\mathcal{S}_{jr}} {p}_{jrs}^t,
~j \in \cN, \: t \in \cT \label{eq:mc-agg}\\
& \mathbf{p}_{jrs} \in \mathcal{P}_{jrs}, \: r\in\mathcal{R}_j, s\in\mathcal{S}_{jr},
~j \in \cN. \label{eq:mc-appl}
\end{align}
\end{subequations}
Linear equality~\eqref{eq:mc-bus} represents the \emph{nodal balance} constraint.
Limits of generator outputs and ramping rates are specified in constraints~\eqref{eq:mc-genlim}
and~\eqref{eq:mc-ramp}. Network line flow constraints are accounted for in~\eqref{eq:mc-flowlim}.
Without loss of generality, the first bus can be set as the reference bus with
zero phase~\eqref{eq:mc-refbus}. Constraint~\eqref{eq:mc-DRAlim} captures the
lower and upper bounds on the energy consumed by the aggregators.
Equality~\eqref{eq:mc-agg} amounts to the \emph{aggregator-users power balance} equation;
finally, \eqref{eq:mc-appl} gives the smart appliance constraints.

A smart appliance example is charging a PHEV battery, which typically amounts to consuming
a prescribed total energy $E_{jrs}$ over a specific horizon from a start
time $T_{jrs}^{\mathrm{st}}$ to a termination time $T_{jrs}^{\mathrm{end}}$.
The consumption must remain within a range between $p_{jrs}^{\mathrm{min}}$
and $p_{jrs}^{\mathrm{max}}$ per period.
With $\mathcal{T}_{E} := \{T_{jrs}^{\mathrm{st}},\ldots,T_{jrs}^{\mathrm{end}}\}$,
set $\mathcal{P}_{jrs}$ takes the form
\begin{align}
\mathcal{P}_{jrs} = \Biggl\{ \mathbf{p}_{jrs} \Biggl| \Biggr.
\sum_{t \in \mathcal{T}_{E}} {p}_{jrs}^t &= E_{jrs};\,
p_{jrs}^t \in [p_{jrs}^{\mathrm{min}}, p_{jrs}^{\mathrm{max}}],\,
\forall~t \in \mathcal{T}_{E}; \nonumber \\
& p_{jrs}^t = 0,~\forall~t \in \cT \setminus \mathcal{T}_{E} \Biggr\}.
\label{eq:phev}
\end{align}
Further examples of $\mathcal{P}_{jrs}$ and $B_{jrs}(\mathbf{p}_{jrs})$ can be found in~\cite{GaGG-TSG13}, where it is argued that $\mathcal{P}_{jrs}$ is a convex set for several appliance types of interest.

Matrices $\mathbf{B}$ and $\mathbf{H}$ are defined as follows.
With $X_{mn}$ denoting the reactance of line $(m,n)$, the bus admittance
matrix $\mathbf{B} \in \mathbb{R}^{N_b \times N_b}$ has elements
\begin{equation*}
\label{eq:Bdef}
[\mathbf{B}]_{mn} = - X_{mn}^{-1}, \:\text{if $m\neq n $};~\quad
[\mathbf{B}]_{mm} = \sum_{n=1}^{N_b}X_{mn}^{-1}
\end{equation*}
where $X_{mn}^{-1}:=0$ if line $(m,n)$ does not exist.
Matrix $\mathbf{H} \in \mathbb{R}^{N_l\times N_b}$ has entries
so that if line $q=1,\ldots,N_l$ connects buses $n$ and $n'$, then
\begin{equation*}
\label{eq:Hdef}
[\mathbf{H}]_{qm} = \begin{cases} X_{nn'}^{-1}, \text{~if $m=n$} \\
 - X_{nn'}^{-1}, \text{~if $m=n'$}  \\
 0, \text{~otherwise}.
\end{cases}
\end{equation*}

Finally, examples detailing the entries of matrices $\mathbf{A}_g$ and $\mathbf{A}_a$
in~\eqref{eq:mc-bus} can be found in~\cite{GaGG-TSG13}.

Problem~\eqref{eq:mc-ALL} can be principally solved at the MO in a central fashion.
However, there are two major challenges when it comes to solving~\eqref{eq:mc-ALL}
with large-scale DR:
i) functions $B_{jrs}(\mathbf{p}_{jrs})$ and sets $\mathcal{P}_{jrs}$ are private,
and cannot be revealed to the MO; ii) including the sheer number of variables
$\mathbf{p}_{jrs}$ would render the overall problem intractable for the MO,
regardless of the privacy issue.
The aggregator plays a critical role in successfully addressing
these two challenges through decomposing the optimization tasks that arises, 
as detailed in the ensuing section.

\section{Decomposition Algorithm}
\label{sec:decomp}

\subsection{Dual Decomposition}
\label{subsec:alg}

Leveraging the dual decomposition technique,
problem~\eqref{eq:mc-ALL} can be decoupled into simpler subproblems tackled by
the MO and the aggregators. Specifically, consider dualizing the
linear coupling constraint~\eqref{eq:mc-agg}
with corresponding Lagrange multiplier $\mu_j^t$.
Upon straightforward re-arrangements, the partial Lagrangian can be written as
\begin{equation}
L(\{\mathbf{p}_j\}_{j=0}^{N_a},\bm{\mu}) =
L_0(\mathbf{p}_0,\bm{\mu}) + \sum_{j=1}^{N_a} L_j(\mathbf{p}_j,\bm{\mu})
\label{eq:LagrDecomp}
\end{equation}
where
\begin{align}
L_0(\mathbf{p}_0,\bm{\mu}) & := \sum_{t=1}^T\left[\sum_{i=1}^{N_g} C_{i}(P_{G_i}^t)
- \sum_{j=1}^{N_a} \mu_{j}^t P_{\mathrm{DRA}_{j}}^t\right] \label{eq:L0}\\
L_j(\mathbf{p}_j,\bm{\mu}) & :=
\sum_{r\in\mathcal{R}_{j}} \sum_{s\in\mathcal{S}_{jr}} \left[
\sum_{t=1}^T \mu_{j}^t  {p}_{jrs}^t -   B_{jrs}(\mathbf{p}_{jrs}) \right].
\label{eq:Lj}
\end{align}

The dual function is thus obtained by minimizing the partial Lagrangian over the
primal variables $\{\mathbf{p}_j\}_{j=0}^{N_a}$ as
\begin{subequations}
\label{eq:DualFn}
\begin{align}
D(\bm\mu) :&=D_0(\bm\mu)+\sum_{j=1}^{N_a} D_j(\bm\mu) \\
&=\min_{\text{s.t. \eqref{eq:mc-bus}--\eqref{eq:mc-DRAlim}}}
L_0(\mathbf{p}_0,\bm{\mu}) + \sum_{j=1}^{N_a} \min_{\text{s.t. \eqref{eq:mc-appl}}} L_j(\mathbf{p}_j,\bm{\mu}).
\end{align}
\end{subequations}

The dual decomposition essentially iterates between two steps:
S1) Lagrangian minimization with respect to $\{p_j\}_{j=0}^{N_a}$ given
the current multipliers, and
S2) multiplier update, using the obtained primal minimizers.
It is clear from~\eqref{eq:LagrDecomp} that the Lagrangian minimization
can be decoupled into $1+N_a$ minimizations,
where one is performed by the MO, and the rest
by the corresponding aggregators.

Specifically, let $k=1,2,\ldots$ index iterations.
Given the multipliers $\bm\mu(k)$, the subproblems at iteration $k$
solved by the MO and each residential end-user are given as follows
\vspace{.2cm}

\begin{subequations}\label{eq:subprobs}
\hspace{-.5cm}
\fbox{
 \addtolength{\linewidth}{-2\fboxsep}%
 \addtolength{\linewidth}{-2\fboxrule}%
 \begin{minipage}{0.98\linewidth}

\begin{align}
\mathbf{p}_0(k)
&= \argmin_{\mathbf{p}_0,~\text{s.t.~\eqref{eq:mc-bus}--\eqref{eq:mc-DRAlim}}}
L_0(\mathbf{p}_0,\bm{\mu}(k)) \label{eq:subMO}  \\
\{\mathbf{p}_{jrs}(k)\}_{s}
&=\argmin_{\{\mathbf{p}_{jrs}\in \mathcal{P}_{jrs}\}_{s}}
\sum_{s\in\mathcal{S}_{jr}} \Big[\sum_{t=1}^T \mu_{j}^t(k) {p}_{jrs}^t
- B_{jrs}(\mathbf{p}_{jrs})\Big]. \label{eq:subUser}
\end{align}
\end{minipage}
}
\end{subequations}

\vspace{.2cm}

Note that subproblem~\eqref{eq:subMO} is a standard DC-OPF while the convex
subproblem~\eqref{eq:subUser} can be handled efficiently by the smart meters.
In fact, with the feasible set in~\eqref{eq:phev} and upon setting
$B_{jrs}(\mathbf{p}_{jrs}) \equiv 0$,~\eqref{eq:subUser} boils down to the
\emph{fractional knapsack} problem, which can be solved in closed form.
To this end, the multipliers $\mu_{j}^t(k)$ needed can be transmitted
to the user's smart meter via the AMI.

With the obtained quantities of $\mathbf{p}_0(k)$,
$\{\mathbf{p}_{jrs}(k)\}_{s}$, and $\{D_j(\bm{\mu}(k))\}_{j=0}^{N_a}$,
the ensuing section develops the approach to updating
the multipliers $\{\mu_j^t\}_{j,t}$ using the so-termed bundle methods.

\subsection{Multiplier Update via Bundle Methods}
\label{subseq:multupd}

The choice of the multiplier update method is crucial, because fewer update
steps imply less communication between the CPM and the aggregators.
A popular method of choice in the context of dual decomposition
is the subgradient method, which is very slow typically.
In this paper, the bundle method with disaggregated cuts is proposed
for the multiplier update. It is better suited to
the problem of interest yielding faster convergence,
because it exploits the special structure of the dual function which can be
written as a sum of separate terms [cf.~\eqref{eq:DualFn}], while it overcomes
the drawbacks of the cutting plane one developed in~\cite{GaGG-TSG13}.
Numerical tests in Section~\ref{sec:num} illustrate
differences in terms of convergence speed.

The following overview of the disaggregated bundle method in a
general form is useful to grasp its role in the present context;
see e.g.,~\cite[Ch.~6]{Bertsekas09} for detailed discussions.
Consider the following separable convex minimization problem
with $n_c$ linear constraints:
\begin{subequations}
\label{eq:genpr-ALL}
\begin{align}
f^*=\min_{\{\mathbf{x}_j \in \mathcal{X}_j\}_{j=0}^{N_a}}
~& \sum_{j=0}^{N_a} f_j(\mathbf{x}_j) \\
\text{s. t.}\quad\,
&\sum_{j=0}^{N_a}\mathbf{A}_j \mathbf{x}_j = \bm{0}.
\label{eq:genpr-coupl}
\end{align}
\end{subequations}
For problem~\eqref{eq:mc-ALL}, constraint~\eqref{eq:genpr-coupl}
corresponds to~\eqref{eq:mc-agg}. Set $\mathcal{X}_0$ captures constraints~\eqref{eq:mc-bus}--\eqref{eq:mc-DRAlim},
while $\mathcal{X}_j$, $j \in \cN$, corresponds to \eqref{eq:mc-appl}.

The dual function $D(\bm\mu) = \sum_{j=0}^{N_a} D_j(\bm\mu)$ can be obtained by
dualizing constraint~\eqref{eq:genpr-coupl} with the multiplier vector $\bm{\mu}$.
Thus, the dual problem is to maximize the dual objective as
\begin{equation}\label{eq:genDualPr}
\max_{\bm\mu\in\mathbb{R}^{n_c}} \sum_{j=0}^{N_a} D_j(\bm\mu) =
\max_{\bm\mu\in\mathbb{R}^{n_c}} \sum_{j=0}^{N_a} \left[\min_{\mathbf{x}_j} \{f_j(\mathbf{x}_j) + \bm\mu'\mathbf{A}_j \mathbf{x}_j\}\right]
\end{equation}
where strong duality holds here due to the polyhedral feasible set~\eqref{eq:genpr-coupl}.

The basic idea of bundle methods (also CPM) is to approximate
the epigraph of a convex (possibly non-smooth) objective function
as the intersection of a number of supporting hyperplanes (also called cuts in this context).
The approximation is gradually refined by generating additional cuts
based on {subgradients} of the objective function.

Specifically, suppose that the method has so far generated
the iterates $\{\bm\mu(\ell)\}_{\ell=1}^k$ after $k$ steps.
Let $\mathbf{x}_j(\ell)$ be the primal minimizer
corresponding to $\bm\mu(\ell)$.
Observe that the vector $\mathbf{g}_j(\ell):= \mathbf{A}_j \mathbf{x}_j(\ell)$
is a subgradient of function $D_j(\bm\mu)$ at point $\bm\mu(\ell)$, and
it thus holds for all $\bm\mu$ such that
\begin{equation}
\label{eq:sgIneq}
D_j(\bm\mu) \leq D_j(\bm\mu(\ell)) + (\bm\mu-\bm\mu(\ell))'\mathbf{g}_j(\ell).
\end{equation}
Clearly, the minimum of the right-hand side of~\eqref{eq:sgIneq}
over $\ell=1,\ldots,k$ is a polyhedral approximation of $D_j(\bm\mu)$,
and is essentially a concave and piecewise linear overestimator of the dual function.

The bundle method with disaggregated cuts generates a sequence
$\{\bm{\mu}(k)\}$ with guaranteed convergence to an optimal solution.
Specifically, the iterate $\bm{\mu}(k+1)$ is obtained by maximizing the
polyhedral approximations of $\{D_j(\bm\mu)\}_{j}$ with a proximal regularization
\begin{subequations}\label{prob:bundle}
\begin{align}
D_{\mathrm{ap}}(\bm{\mu}(k+1)) := \max_{\bm\mu,\{v_j\}_{j=1}^{N_a}}
&\sum_{j=0}^{N_a} v_j - \frac{\rho(k)}{2}\|\bm{\mu}-\check{\bm{\mu}}(k)\|_2^2 \label{eq:bundleObj} \\
\text{s. t.}~\quad & v_j \leq D_j(\bm\mu(\ell)) + (\bm\mu-\bm\mu(\ell))'\mathbf{g}_j(\ell), \notag\\
& \,\, j=0,\ldots,N_a, \ell=1,\ldots,k
\end{align}
\end{subequations}
where the proximity weight $\rho(k)>0$ is to control stability of the iterates;
and the proximal center $\check{\bm{\mu}}(k)$ is updated according to a query
for ascent
\begin{align*}
\check{\bm{\mu}}(k+1)  =
\left\{\begin{array}{cc}
\bm{\mu}(k+1),  &\mbox{if}~D(\bm{\mu}(k+1))-D(\check{\bm{\mu}}(k)) \ge \beta \eta(k) \\
\check{\bm{\mu}}(k),  &\mbox{otherwise}
\end{array}\right.
\end{align*}
where $\eta(k) = D_{\mathrm{ap}}(\bm{\mu}(k+1))-D(\check{\bm{\mu}}(k))$,
and $\beta \in (0,1)$.
Finally, the bundle algorithm can be terminated when
$\eta(k)<\epsilon$ holds for
a prescribed tolerance $\epsilon$ (cf.~\cite[Ch.~6]{Bertsekas09}).

\begin{remark}\textit{(Bundle methods versus CPM)}.
When $\rho(k) \equiv 0$, problem~\eqref{prob:bundle} boils down to
the CPM with disaggregated cuts for solving the dual, which is however known to be unstable and
converges slowly on some practical instances~\cite{Lemarechal93}.
The proximal regularization in the bundle methods is thus introduced
to improve stability of the iterates, while the \emph{smart} prox-center
updating rule enhances further the convergence speed compared with
the proximal CPM. A further limitation of CPM is that a compact set containing the optimal solution
has to be included, as is the case with  $\bm{\mu}\in [\bm{\mu}^{\min},\bm{\mu}^{\max}]$ in~\cite{GaGG-TSG13}.
The CPM convergence performance depends on the choice of this set, while there is no such issue for the bundle methods.
Note further that the dual problem of~\eqref{prob:bundle}
is a quadratic program (QP) over a probability simplex.
Such a special structure can be exploited by off-the-shelf
QP solvers, and hence it is efficiently solvable.
As a result, solving~\eqref{prob:bundle} does not
require much more computational work than solving a linear program (LP),
which is the case for the CPM.
Finally, it is worth stressing that the disaggregated bundle method
takes advantage of the \emph{separability} of~\eqref{eq:genpr-ALL}.
In a nutshell, offering state-of-the-art algorithms for solving
non-smooth convex programs, the stable and fast convergent bundle methods
are well motivated here for clearing the market distributedly.
\end{remark}

Specifically, applying the disaggregated bundle method to
problem~\eqref{eq:genDualPr} at hand, the multiplier update
at iteration $k$ amounts to solving the following problem:
 \vspace{.2cm}

\begin{subequations}\label{eq:bundleMC}
\hspace{-.5cm}
\fbox{
 \addtolength{\linewidth}{-2\fboxsep}%
 \addtolength{\linewidth}{-2\fboxrule}%
 \begin{minipage}{0.98\linewidth}
\begin{align}
\max_{\{\mu_j^t, v_j\}_{j,t}}
&\sum_{j=0}^{N_a} v_j - \frac{\rho(k)}{2}
\sum_{j=1}^{N_a}\sum_{t=1}^T(\mu_j^t-\check{\mu}_j^t)^2 \label{eq:bundleObjMC} \\
\text{s. t.}~\quad
&v_0 \leq D_0(\bm{\mu}(\ell))
- \sum_{j=1}^{N_a}\sum_{t=1}^TP_{\mathrm{DRA}_j}^t(\ell)[\mu_j^t-\mu_j^t(\ell)] \label{eq:v0cstt}
\nonumber \\  &\hspace{3.7cm}  \ell=1,\ldots,k  \\
& v_j \leq D_j(\bm{\mu}(\ell))
+ \sum_{t=1}^T\sum_{r,s}p_{jrs}^t(\ell)[\mu_j^t-\mu_j^t(\ell)]\label{eq:vjcstt}
\nonumber \\
& \hspace{2.5cm} j \in \cN,~\ell=1,\ldots,k.
\end{align}
\end{minipage}
}
\end{subequations}

\vspace{.2cm}

\begin{figure}[t]
\centering
\includegraphics[scale=0.32]{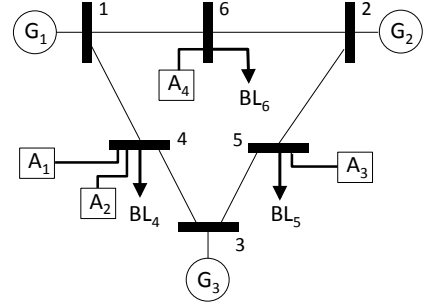}
\caption{Power system example featuring 6 buses, 3 generators, 4 aggregators, and base loads at three of the buses.}
\label{fig:wecc}
\vspace{-0.4cm}
\end{figure}

Problem~\eqref{eq:bundleMC} that yields the updated multipliers $\bm{\mu}(k+1)$
and the approximate dual value $D_{\mathrm{ap}}(\bm{\mu}(k+1))$
can be solved at the MO. To this end, the quantities
$\{D_j(\bm{\mu}(k)), \sum_{r,s}{p}_{jrs}^t(k)\}_{j}$
are needed from each aggregator per iteration $k$
as the problem input. Note that
$D_j(\bm{\mu}(k)) := \sum_{r\in \mathcal{R}_j}D_{jr}(\bm{\mu}(k))$,
where $D_{jr}(\bm{\mu}(k))$ is the optimal value of problem~\eqref{eq:subUser}.
Thus, it is clear that all these required quantities can be formed at
the aggregator level as summations over all end-users, and then
transmitted to the MO.
The highlight here is that the proposed decomposition scheme
respects \emph{user privacy}, since $B_{jrs}(\mathbf{p}_{jrs})$
and $\mathcal{P}_{jrs}$ are never revealed.


\section{Numerical tests}
\label{sec:num}

In this section, simulation results are presented to verify the
merits of the disaggregated bundle method.
The power system tested for market clearing and large-scale DR is
illustrated in Fig.~\ref{fig:wecc}, where each of the 4 aggregators
serves 1,000 residential end-users. The scheduling horizon starts
from 1am until 12am, for a total of 24 hours.

Time-invariant generation cost functions are set to be quadratic as
$C_i(P_{G_i}^t)=a_i (P_{G_i}^t)^2 + b_iP_{G_i}^t$ for all $i$ and $t$.
Each end-user has a PHEV to charge overnight. All detailed
parameters of the generators and loads are listed in Tables~\ref{T:gen} and~\ref{T:load}.
The utility functions $\{B_{jrs}(\cdot)\}$ are set to be zero for simplicity.
The upper bound on each aggregator's consumption is
$P_{\mathrm{DRA}_j}^{\max}=50$ MW while $\mathbf{p}_{\mathrm{BL}}^t = 5$ MW.
At a base of 100 MVA, the values of the network reactances are
$\{X_{16},X_{62},X_{25},X_{53},X_{34},X_{41}\}=\{0.2, 0.3, 0.25, 0.1, 0.3, 0.4\}$ p.u.
Finally, no flow limits are imposed across the network.
The resulting optimization problems~\eqref{eq:subMO} and~\eqref{eq:bundleMC}
are modeled via~\texttt{YALMIP}~\cite{YALMIP}, and solved by \texttt{Gurobi}~\cite{gurobi}.

\begin{table}[t]
\centering
\caption{Generator parameters.
The units of $a_i$ and $b_i$ are \$/(MWh)$^{2}$ and \$/MWh, respectively.
the rest are in MW.}
\label{T:gen}
\begin{tabular}{c|c|c|c|c|c}
Gen. &$a_i$ &$b_i$ &$P_{G_i}^{\max}$ & $P_{G_i}^{\min}$ & $\mathsf{R}_i^{\mathrm{up,down}}$ \\
\hline
1 &0.3  & 3  & 60 & 2.4 & 50  \\
2 &0.15  & 20 & 50 & 0   & 35  \\
3 &0.2   & 50 & 50 & 0   & 40  \\
\hline
\end{tabular}
\end{table}

\begin{table}[t]
\renewcommand{\arraystretch}{1.1}
\centering
\caption{Parameters of residential appliances.
All listed hours are the ending ones;
w.p. means with probability.}
\label{T:load}
\begin{tabular}{ c || c }
\hline
$E_{\mathrm{PHEV}}$  (kWh)       & Uniform on \{10, 11, 12\} \\
$p_{\mathrm{PHEV}}^{\max}$ (kWh) & Uniform on \{2.1, 2.3, 2.5\}  \\
$p_{\mathrm{PHEV}}^{\min}$ (kWh) &  0   \\
$T_{jr1}^{\mathrm{st}}$          & 1am \\
$T_{jr1}^{\mathrm{end}}$         & 6am w.p. 70\%, 7am w.p. 30\% \\
\hline
\end{tabular}
\end{table}

\begin{figure}[t]
\centering
\includegraphics[scale=.47]{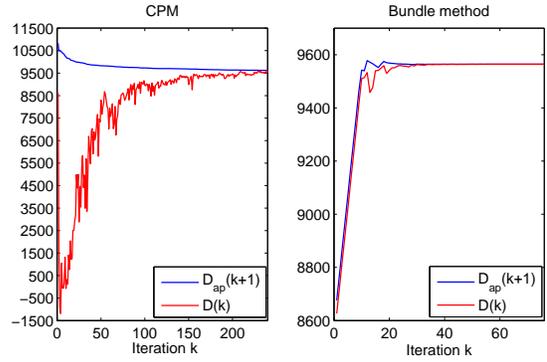}
\vspace{-0.2cm}
\caption{Convergence of the objective values of
the dual and the approximated one
(denoted as $D(k)$ and $D_{\mathrm{ap}}(k+1)$ in the caption).}
\label{fig:objConv}
\end{figure}

\begin{figure}[t]
\centering
\includegraphics[scale=.47]{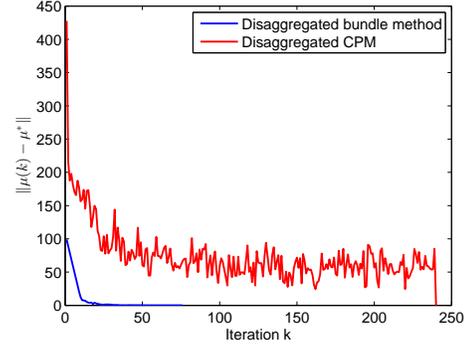}
\vspace{-0.2cm}
\caption{Convergence of the Lagrange multipliers.}
\label{fig:argConv}
\vspace{-0.4cm}
\end{figure}

Figs.~\ref{fig:objConv} and~\ref{fig:argConv} illustrate the
convergence performance of the proposed disaggregated
bundle method vis-\`{a}-vis the disaggregated CPM.
The pertinent parameters are set as $\epsilon=10^{-3}$,
$\rho(k)\equiv0$, $\beta = 0.5$,
and $\bm\mu^{\max,\min}=\mathbf{\pm 50}$ (cf.~\cite{GaGG-TSG13}).
Fig.~\ref{fig:objConv} depicts the evolution of the objective values of the
dual $D(\bm{\mu}(k))$ and the approximate dual $D_{\mathrm{ap}}(\bm{\mu}(k+1))$.
It is clearly seen that the bundle method converges much faster (more than three times)
than its CPM counterpart. Note that due to the effect of the proximal penalty (cf.~\eqref{eq:bundleObj}),
quantity $D_{\mathrm{ap}}(\bm{\mu}(k+1))$ for the bundle
may not always serve as an upper bound of $f^{*}$ as the one for the CPM.
Finally, convergence of the Lagrange multiplier sequence $\bm{\mu}(k)$ is
shown in Fig.~\ref{fig:argConv}, which also corroborates the merit of the bundle method for
its faster parameter convergence over the CPM. It is interesting to
observe that the distance-to-optimal curve of the bundle method
is quite smooth compared with the CPM one. This again illustrates the effect of the
proximal regulation penalizing large deviations.

\section{Conclusions and Future Directions}
\label{sec:concl}

In this work, a fast convergent and scalable distributed solver is developed
for market clearing with large-scale residential DR.
Leveraging the dual decomposition technique, only the aggregator-users balance
constraint is dualized in order to separate problems for the MO and each aggregator,
while respecting end-user privacy concerns. Simulated tests highlight the merits of
the proposed approach for multiplier updates based on the disaggregated bundle method.

A number of interesting research directions open up, including the incorporation of 
load and renewable energy production uncertainty, the issue of primal recovery, 
as well as cut aggregation techniques for further computational speed up.


\bibliographystyle{IEEEtran}
\bibliography{MCDR_biblio}

\end{document}